\documentclass{amsart}

\usepackage{amssymb,amsmath,amsthm,latexsym,booktabs}

\theoremstyle{definition}
\newtheorem{definition}{Definition}[section]

\theoremstyle{plain}
\newtheorem{lemma}[definition]{Lemma}
\newtheorem{theorem}[definition]{Theorem}

\begin{document}

\title{Higher identities for the ternary commutator}

\author{Murray R. Bremner}

\address{Department of Mathematics and Statistics, University of Saskatchewan, Canada}

\email{bremner@math.usask.ca}

\author{Luiz A. Peresi}

\address{Department of Mathematics, University of S\~ao Paulo, Brazil}

\email{peresi@ime.usp.br}

\date{25 July 2012}



\begin{abstract}
We use computer algebra to study polynomial identities for the trilinear operation
$[a,b,c] = abc - acb - bac + bca + cab - cba$ in the free associative algebra.
It is known that $[a,b,c]$ satisfies the alternating property in degree 3,
no new identities in degree 5, a multilinear identity in degree 7 which alternates
in 6 arguments, and no new identities in degree 9.
We use the representation theory of the symmetric group to demonstrate the existence
of new identities in degree 11.
The only irreducible representations of dimension $< 400$ with new identities correspond
to partitions $2^5 1$ and $2^4 1^3$ and have dimensions 132 and 165.
We construct an explicit new multilinear identity for partition $2^5 1$
and we demonstrate the existence of a new non-multilinear identity
in which the underlying variables are permutations of $a^2b^2c^2d^2e^2f$.
\end{abstract}

\maketitle


\section{Introduction}

The theory of multioperator algebras ($\Omega$-algebras), by which is meant vector spaces
with multilinear operations, was first studied systematically by the school of Kurosh in Moscow; see \cite{BB,K}.
In particular, a natural generalization of the Lie bracket to the $n$-ary setting is the
alternating $n$-ary sum ($n$-commutator):
  \[
  [ a_1, \dots, a_n ]
  =
  \sum_{\sigma \in S_n}
  \epsilon(\sigma) \,
  a_{\sigma(1)} \cdots a_{\sigma(n)},
  \]
where $\epsilon(\sigma)$ is the sign of the permutation $\sigma$.
This operation provides unexpected algebraic structures on vector fields \cite{D1,D2}, and plays
an essential role in the construction of universal enveloping algebras of Filippov algebras
($n$-Lie algebras) \cite{EB}.
For many other applications, especially to theoretical physics, see the recent survey of
$n$-ary analogues of Lie algebras \cite{DI}.

For $n = 3$, the alternating ternary sum (ternary commutator) has the form
  \[
  [a,b,c] = abc - acb - bac + bca + cab - cba.
  \]
The first explicit polynomial identity which is satisfied by this operation in the free associative algebra,
but which does not follow from the alternating property in degree 3, was found in 1998; see \cite{B2}.
This identity has degree 7:
  \[
  \sum_{\sigma \in S_6}
  \epsilon(\sigma) \,
  \big( \,
  [[[b^\sigma,c^\sigma,d^\sigma],a,e^\sigma],f^\sigma,g^\sigma]
  +
  [[a,b^\sigma,c^\sigma],[d^\sigma,e^\sigma,f^\sigma],g^\sigma]
  \, \big)
  \equiv
  0.
  \]
Two years later, it was shown that there are no new identities in degree 9; see \cite{BH}.
Ten years later, the identity in degree 7 was rediscovered \cite{DFNW},
and was generalized to all odd $n$ \cite{CJM};
the situation is much simpler for even $n$ \cite{G,HW,R}.
For identities relating the ordinary and ternary commutators, see \cite{FN}.
For the partially alternating ternary sum in an associative dialgebra (Loday algebra), see \cite{BSO}.

In this paper, we use computer algebra to show that further new polynomial identities for
the ternary commutator exist in degree 11.
We construct an explicit new multilinear identity for partition $2^5 1$
and we demonstrate the existence of a new non-multilinear identity
in which the underlying variables are permutations of $a^2b^2c^2d^2e^2f$.

Owing to the large size of the matrices involved in our computations, we used modular arithmetic
to save memory.  By choosing a suitable modulus, we found it easy to perform ``rational reconstruction''
of the correct results in characteristic 0 from the results obtained in characteristic $p$.  Underlying all
of these computations is the structure theory of the group algebra of the symmetric group $S_n$, which is
semisimple both in characteristic 0 and in characteristic $p > n$.  For further information, see
\cite[Lemma 8]{BP1} and \cite[\S 5.5]{BP2}.


\section{Preliminaries} \label{preliminaries}

For an alternating trilinear operation, every monomial in degree 11 can be written in terms of
one of the following eight association types (placements of brackets);
we display these types with the identity permutation of the arguments:
  \begin{equation}
  \label{types} \tag{$\ast$}
  \left\{
  \begin{array}{ll}
  1\colon \; [[[[[a,b,c],d,e],f,g],h,i],j,k] &\quad
  2\colon \; [[[[a,b,c],[d,e,f],g],h,i],j,k] \\[3pt]
  3\colon \; [[[[a,b,c],d,e],[f,g,h],i],j,k] &\quad
  4\colon \; [[[a,b,c],[d,e,f],[g,h,i]],j,k] \\[3pt]
  5\colon \; [[[[a,b,c],d,e],f,g],[h,i,j],k] &\quad
  6\colon \; [[[a,b,c],[d,e,f],g],[h,i,j],k] \\[3pt]
  7\colon \; [[[a,b,c],d,e],[[f,g,h],i,j],k] &\quad
  8\colon \; [[[a,b,c],d,e],[f,g,h],[i,j,k]]
  \end{array}
  \right.
  \end{equation}
The number of multilinear monomials in each type can be easily calculated using the
alternating property of $[a,b,c]$; the total is
  \[
                      \frac{11!}{6^1 2^4} +
  \frac{1}{2} {\cdot} \frac{11!}{6^2 2^2} +
                      \frac{11!}{6^2 2^2} +
  \frac{1}{6} {\cdot} \frac{11!}{6^3 2} +
                      \frac{11!}{6^2 2^2} +
  \frac{1}{2} {\cdot} \frac{11!}{6^3} +
  \frac{1}{2} {\cdot} \frac{11!}{6^2 2^2} +
  \frac{1}{2} {\cdot} \frac{11!}{6^3 2}
  =
  1401400.
  \]
Since this number is so large, we cannot process all the monomials at once, so we
use the representation theory of the symmetric group $S_{11}$ to decompose the problem into
a sequence of smaller problems, each corresponding to an irreducible representation.
(For a detailed discussion of this approach, see \cite[\S4]{BP1} or \cite[\S5]{BP2}.)

Using representation theory requires that we enumerate the symmetries of the association types \eqref{types};
each symmetry is a two-term identity expressing the fact that the value of a monomial
changes sign after a transposition of two factors.
Since the symmetric group $S_3$ is generated by the transpositions $(12)$ and $(23)$,
every symmetry is a consequence of the 43 symmetries corresponding to the monomials $\pi$
in Table \ref{degree11symmetries}.
In that table, $\pi$ represents the identity
$\iota + \pi \equiv 0$, where $\iota$ represents the monomial with the identity permutation
of the variables in the same association type.
These symmetries are the consequences in degree 11 of the alternating properties
$[a,b,c] + [b,a,c] \equiv 0$ and $[a,b,c] + [a,c,b] \equiv 0$ in degree 3.

  \begin{table}
  \[
  \begin{array}{cc}
  1\left\{
  \begin{array}{l}
  {[[[[[b,a,c],d,e],f,g],h,i],j,k]} \\
  {[[[[[a,c,b],d,e],f,g],h,i],j,k]} \\
  {[[[[[a,b,c],e,d],f,g],h,i],j,k]} \\
  {[[[[[a,b,c],d,e],g,f],h,i],j,k]} \\
  {[[[[[a,b,c],d,e],f,g],i,h],j,k]} \\
  {[[[[[a,b,c],d,e],f,g],h,i],k,j]}
  \end{array}
  \right.
  &
  2\left\{
  \begin{array}{l}
  {[[[[b,a,c],[d,e,f],g],h,i],j,k]} \\
  {[[[[a,c,b],[d,e,f],g],h,i],j,k]} \\
  {[[[[d,e,f],[a,b,c],g],h,i],j,k]} \\
  {[[[[a,b,c],[d,e,f],g],i,h],j,k]} \\
  {[[[[a,b,c],[d,e,f],g],h,i],k,j]}
  \end{array}
  \right.
  \\[36pt]
  3\left\{
  \begin{array}{l}
  {[[[[b,a,c],d,e],[f,g,h],i],j,k]} \\
  {[[[[a,c,b],d,e],[f,g,h],i],j,k]} \\
  {[[[[a,b,c],e,d],[f,g,h],i],j,k]} \\
  {[[[[a,b,c],d,e],[g,f,h],i],j,k]} \\
  {[[[[a,b,c],d,e],[f,h,g],i],j,k]} \\
  {[[[[a,b,c],d,e],[f,g,h],i],k,j]}
  \end{array}
  \right.
  &
  4\left\{
  \begin{array}{l}
  {[[[b,a,c],[d,e,f],[g,h,i]],j,k]} \\
  {[[[a,c,b],[d,e,f],[g,h,i]],j,k]} \\
  {[[[d,e,f],[a,b,c],[g,h,i]],j,k]} \\
  {[[[a,b,c],[g,h,i],[d,e,f]],j,k]} \\
  {[[[a,b,c],[d,e,f],[g,h,i]],k,j]}
  \end{array}
  \right.
  \\[36pt]
  5\left\{
  \begin{array}{l}
  {[[[[b,a,c],d,e],f,g],[h,i,j],k]} \\
  {[[[[a,c,b],d,e],f,g],[h,i,j],k]} \\
  {[[[[a,b,c],e,d],f,g],[h,i,j],k]} \\
  {[[[[a,b,c],d,e],g,f],[h,i,j],k]} \\
  {[[[[a,b,c],d,e],f,g],[i,h,j],k]} \\
  {[[[[a,b,c],d,e],f,g],[h,j,i],k]}
  \end{array}
  \right.
  &
  6\left\{
  \begin{array}{l}
  {[[[b,a,c],[d,e,f],g],[h,i,j],k]} \\
  {[[[a,c,b],[d,e,f],g],[h,i,j],k]} \\
  {[[[d,e,f],[a,b,c],g],[h,i,j],k]} \\
  {[[[a,b,c],[d,e,f],g],[i,h,j],k]} \\
  {[[[a,b,c],[d,e,f],g],[h,j,i],k]}
  \end{array}
  \right.
  \\[36pt]
  7\left\{
  \begin{array}{l}
  {[[[b,a,c],d,e],[[f,g,h],i,j],k]} \\
  {[[[a,c,b],d,e],[[f,g,h],i,j],k]} \\
  {[[[a,b,c],e,d],[[f,g,h],i,j],k]} \\
  {[[[f,g,h],i,j],[[a,b,c],d,e],k]}
  \end{array}
  \right.
  &
  8\left\{
  \begin{array}{l}
  {[[[b,a,c],d,e],[f,g,h],[i,j,k]]} \\
  {[[[a,c,b],d,e],[f,g,h],[i,j,k]]} \\
  {[[[a,b,c],e,d],[f,g,h],[i,j,k]]} \\
  {[[[a,b,c],d,e],[g,f,h],[i,j,k]]} \\
  {[[[a,b,c],d,e],[f,h,g],[i,j,k]]} \\
  {[[[a,b,c],d,e],[i,j,k],[f,g,h]]}
  \end{array}
  \right.
  \end{array}
  \]
  \caption{The symmetries of the association types in degree 11}
  \label{degree11symmetries}
  \end{table}

We also need to determine the consequences in degree 11 of the known polynomial identity in degree 7;
see \cite{B2,CJM,DFNW}.
We write this identity symbolically as $I(a,b,c,d,e,f,g) \equiv 0$, where
  \begin{align*}
  &
  I(a,b,c,d,e,f,g)
  =
  \\
  &\quad
  \sum_{\sigma \in S_6} \epsilon(\sigma) \,
  \big( \,
  [ [ [ b^\sigma, c^\sigma, d^\sigma ], a, e^\sigma ], f^\sigma, g^\sigma ]
  +
  [ [ a, b^\sigma, c^\sigma ], [ d^\sigma, e^\sigma, f^\sigma ], g^\sigma ]
  \, \big).
  \end{align*}
We collect similar terms in this identity using the alternating property of the ternary commutator,
and see that the total number of distinct terms is
$\binom{6}{3,1,2} + \binom{6}{2,3,1} = 60 + 60 = 120$.
From the alternating property of $[a,b,c]$, and the alternating property of $I(a,b,c,d,e,f,g)$ in the arguments $b, \dots, g$,
it follows that every consequence of $I(a,b,c,d,e,f,g) \equiv 0$ in degree 9 is a linear combination of permutations of
three identities, the first two obtained by substituting a triple for a variable, and the third obtained by embedding
the identity in a triple:
  \allowdisplaybreaks
  \begin{alignat*}{2}
  &I([a,h,i],b,c,d,e,f,g) \equiv 0, &\qquad
  &I(a,[b,h,i],c,d,e,f,g) \equiv 0, \\
  &[I(a,b,c,d,e,f,g),h,i] \equiv 0.
  \end{alignat*}
Similarly, every consequence of these three identities in degree 11 is a linear combination
of permutations of the eight identities in Table \ref{degree11liftings}.
(It is a coincidence that the number of association types is equal to the number of consequences of
$I(a,b,c,d,e,f,g)$.)
We call these consequences the ``liftings'' of $I(a,b,c,d,e,f,g)$ to degree 11.
We summarize this discussion in the following lemma.

  \begin{table}
  \begin{alignat*}{2}
  1\colon &\; I([[a,j,k],h,i],b,c,d,e,f,g) \equiv 0, &\qquad
  2\colon &\; I([a,h,i],[b,j,k],c,d,e,f,g) \equiv 0, \\
  3\colon &\; [I([a,h,i],b,c,d,e,f,g),j,k] \equiv 0, &\qquad
  4\colon &\; I(a,[[b,j,k],h,i],c,d,e,f,g) \equiv 0, \\
  5\colon &\; I(a,[b,h,i],[c,j,k],d,e,f,g) \equiv 0, &\qquad
  6\colon &\; [I(a,[b,h,i],c,d,e,f,g),j,k] \equiv 0, \\
  7\colon &\; [I(a,b,c,d,e,f,g),[h,j,k],i] \equiv 0, &\qquad
  8\colon &\; [[I(a,b,c,d,e,f,g),h,i],j,k] \equiv 0.
  \end{alignat*}
  \caption{The consequences of $I(a,b,c,d,e,f,g) \equiv 0$ in degree 11}
  \label{degree11liftings}
  \end{table}

\begin{lemma}
Every polynomial identity in degree 11 satisfied by the ternary commutator,
which is a consequence of identities of lower degree, is a linear combination of
permutations of the identities in Tables \ref{degree11symmetries} and \ref{degree11liftings}.
\end{lemma}


\section{New identities in degree 11} \label{newidentity}

Let $\lambda$ be a partition of 11 with $k$ parts; we write
  \[
  \lambda = ( \lambda_1, \dots, \lambda_k ), \qquad
  \lambda_1 \ge \cdots \ge \lambda_k, \qquad
  \lambda_1 + \cdots + \lambda_k = 11.
  \]
Let $d_\lambda$ be the dimension of the corresponding irreducible representation of $S_{11}$.
For any permutation $\pi \in S_{11}$ the $d_\lambda \times d_\lambda$ matrix $R_\pi^\lambda$
representing $\pi$ in the natural representation can be computed using the methods of \cite{Clifton};
see also \cite[\S 5]{BP2}.
Taking linear combinations gives the matrix representing any element of the group algebra $\mathbb{Q} S_{11}$
over the rational field $\mathbb{Q}$.
This provides an algorithm for explicit computation of the isomorphism $\phi$ from $\mathbb{Q} S_{11}$
to its Wedderburn decomposition, by which we mean the direct sum over all partitions $\lambda$ of
matrix algebras of size $d_\lambda \times d_\lambda$.
Any multilinear polynomial $P$ of degree 11 in the ternary commutator can be expressed as a
sum of eight elements of $\mathbb{Q} S_{11}$, one term for each association type \eqref{types}.
For each partition $\lambda$, the projection of $P$ to the corresponding component of the Wedderburn
decomposition consists of an ordered list of eight $d_\lambda \times d_\lambda$ matrices.
The partitions $\lambda$ for which $d_\lambda < 400$ are given in Table \ref{degree11ranks}.

  \begin{table}
  \[
  \begin{array}{rrlrrrrl}
  \# &\quad d_\lambda &\quad \lambda &\quad \text{sym} &\quad \text{sym+lif} &\quad \text{all} &\quad \text{new} \\
  \midrule
   1  &\quad    1  &\quad  11          &\quad    8 &\quad    8 &\quad    8 &\quad 0 &\quad \\
   2  &\quad   10  &\quad  10, 1       &\quad   80 &\quad   80 &\quad   80 &\quad 0 &\quad \\
   3  &\quad   44  &\quad  9, 2        &\quad  352 &\quad  352 &\quad  352 &\quad 0 &\quad \\
   4  &\quad   45  &\quad  9, 1^2      &\quad  360 &\quad  360 &\quad  360 &\quad 0 &\quad \\
   5  &\quad  110  &\quad  8, 3        &\quad  880 &\quad  880 &\quad  880 &\quad 0 &\quad \\
   6  &\quad  231  &\quad  8, 2, 1     &\quad 1848 &\quad 1848 &\quad 1848 &\quad 0 &\quad \\
   7  &\quad  120  &\quad  8, 1^3      &\quad  960 &\quad  960 &\quad  960 &\quad 0 &\quad \\
   8  &\quad  165  &\quad  7, 4        &\quad 1320 &\quad 1320 &\quad 1320 &\quad 0 &\quad \\
  10  &\quad  385  &\quad  7, 2^2      &\quad 3080 &\quad 3080 &\quad 3080 &\quad 0 &\quad \\
  12  &\quad  210  &\quad  7, 1^4      &\quad 1680 &\quad 1680 &\quad 1680 &\quad 0 &\quad \\
  13  &\quad  132  &\quad  6, 5        &\quad 1056 &\quad 1056 &\quad 1056 &\quad 0 &\quad \\
  19  &\quad  252  &\quad  6, 1^5      &\quad 2016 &\quad 2016 &\quad 2016 &\quad 0 &\quad \\
  20  &\quad  330  &\quad  5^2, 1      &\quad 2639 &\quad 2639 &\quad 2639 &\quad 0 &\quad \\
  29  &\quad  210  &\quad  5, 1^6      &\quad 1676 &\quad 1676 &\quad 1676 &\quad 0 &\quad \\
  40  &\quad  120  &\quad  4, 1^7      &\quad  944 &\quad  948 &\quad  948 &\quad 0 &\quad \\
  45  &\quad  385  &\quad  3^2, 1^5    &\quad 3005 &\quad 3020 &\quad 3020 &\quad 0 &\quad \\
  46  &\quad  330  &\quad  3, 2^4      &\quad 2639 &\quad 2639 &\quad 2639 &\quad 0 &\quad \\
  49  &\quad  231  &\quad  3, 2, 1^6   &\quad 1764 &\quad 1795 &\quad 1795 &\quad 0 &\quad \\
  50  &\quad   45  &\quad  3, 1^8      &\quad  333 &\quad  349 &\quad  349 &\quad 0 &\quad \\
  51  &\quad  132  &\quad  2^5, 1      &\quad 1006 &\quad 1020 &\quad 1021 &\quad 1 &\quad \leftarrow \\
  52  &\quad  165  &\quad  2^4, 1^3    &\quad 1242 &\quad 1269 &\quad 1270 &\quad 1 &\quad \leftarrow \\
  53  &\quad  110  &\quad  2^3, 1^5    &\quad  807 &\quad  842 &\quad  842 &\quad 0 &\quad \\
  54  &\quad   44  &\quad  2^2, 1^7    &\quad  302 &\quad  333 &\quad  333 &\quad 0 &\quad \\
  55  &\quad   10  &\quad  2, 1^9      &\quad   57 &\quad   76 &\quad   76 &\quad 0 &\quad \\
  56  &\quad    1  &\quad  1^{11}      &\quad    0 &\quad    7 &\quad    7 &\quad 0 &\quad \\
  \midrule
  \end{array}
  \]
  \caption{Representations of $S_{11}$ with dimension $< 400$}
  \label{degree11ranks}
  \end{table}

We apply this discussion to the symmetries of the association types in Table \ref{degree11symmetries}
and the consequences of $I(a,b,c,d,e,f,g)$ in Table \ref{degree11liftings}.
We first construct a $43 d_\lambda \times 8 d_\lambda$ matrix consisting of $d_\lambda \times d_\lambda$
blocks: for $i = 1, \dots, 43$ and $j = 1, \dots, 8$, block $(i,j)$ contains
the representation matrix for the terms of symmetry $i$ in association type $j$.
Since each symmetry has the form $\iota + \pi$ in one association type,
for each $i$ there will be one nonzero block containing the matrix $I + R_\pi^\lambda$.
We compute the row canonical form of this matrix; for each $\lambda$, the nonzero rows
form a basis for the space of identities in degree 11 which are consequences of
the alternating property of the ternary commutator.
The rank $s_\lambda$ of this matrix is given in column ``sym'' of Table \ref{degree11ranks}.

We next construct a $51 d_\lambda \times 8 d_\lambda$ matrix consisting of $d_\lambda \times d_\lambda$ blocks;
the first 43 rows of blocks are the same as in the preceding matrix.
For $i = 1, \dots, 8$ and $j = 1, \dots, 8$, block $(43{+}i,j)$ contains
the representation matrix for the terms in association type $j$ of the $i$-th consequence of $I(a,b,c,d,e,f,g)$.
We compute the row canonical form of this matrix; for each $\lambda$, the nonzero rows
form a basis for the space of identities in degree 11 which are consequences of all the identities of lower degree.
The rank $sl_\lambda$ of this matrix is given in column ``sym+lif'' of Table \ref{degree11ranks}.
In the next section we use the name $\texttt{oldmat}(\lambda)$ for the $sl_\lambda \times 8 d_\lambda$
matrix in row canonical form containing the nonzero rows.

Finally, we construct a matrix of size $8 d_\lambda \times 9 d_\lambda$ and use it to find all
the identities satisfied by the ternary commutator in degree 11.
The first column of $d_\lambda \times d_\lambda$ blocks corresponds to the associative multilinear polynomials,
which we identify with the group algebra $\mathbb{Q} S_{11}$.
The remaining columns correspond to the eight association types \eqref{types}.
For $i = 1, \dots, 8$ we put the identity matrix in block $(i,i{+}1)$;
in block $(i,1)$ we put the representation matrix for the expansion of association type $i$ (with the identity
permutation of the variables) in the free associative algebra.
By the expansion of an association type, we mean the associative polynomial obtained by replacing
each occurrence of $[a,b,c]$ by the alternating ternary sum of its arguments;
thus each expansion is a sum of $6^5 = 7776$ terms with coefficients $\pm 1$.
In the resulting matrix, the $8 d_\lambda \times 8 d_\lambda$ submatrix obtained by deleting the first column of blocks
is the identity matrix; hence the matrix has rank $8 d_\lambda$.
We compute the row canonical form of this matrix, and delete the rows whose leading 1s occur
within the first $d_\lambda$ columns.
From the remaining matrix, we delete the first $d_\lambda$ columns, all of whose entries are 0.
The result is a matrix of size $a_\lambda \times 8 d_\lambda$, with rank $a_\lambda$ for some $a_\lambda \ge 0$;
this number is given in column ``all'' of Table \ref{degree11ranks}.
For each $\lambda$, the (nonzero) rows of this matrix
provide a basis for the space of all identities in degree 11 satisfied by the ternary commutator.
In the next section we call this matrix $\texttt{allmat}(\lambda)$.

It is clear that $a_\lambda \ge sl_\lambda$ for every $\lambda$: the space of identities
which are consequences of identities of lower degree is a subspace of the space of all identities.
If $a_\lambda = sl_\lambda$ for some $\lambda$ then there are no new identities for partition $\lambda$.
In this case we also need to verify that the two matrices are exactly the same:
the first matrix, $\texttt{oldmat}(\lambda)$ of size $sl_\lambda \times 8 d_\lambda$, containing the symmetries of the association types
and the consequences of $I(a,b,c,d,e,f,g)$;
and the second matrix, $\texttt{allmat}(\lambda)$ of size $a_\lambda \times 8 d_\lambda$, containing all the identities satisfied by
the ternary commutator.
If $a_\lambda > sl_\lambda$ for some $\lambda$ then there exist new identities in degree 11
for the representation of $S_{11}$ corresponding to $\lambda$.
The difference $a_\lambda - sl_\lambda$ is given in column ``new'' of Table \ref{degree11ranks}.

Owing to the large size of many of the irreducible representations of $S_{11}$,
and the time required to compute the representation matrices $R^\lambda_\pi$,
we were able to complete these computations only for the 25 partitions in Table \ref{degree11ranks},
corresponding to the representations with dimensions $< 400$,
slightly less than half of the total of 56 representations.
We found two representations which have new identities: number 51 (dimension 132, partition $2^5 1$) and
number 52 (dimension 165, partition $2^4 1^3$).
We summarize this discussion in the following theorem.

\begin{theorem} \label{newtheorem}
New identities in degree 11 for the ternary commutator exist for partitions $2^51$ and $2^4 1^3$,
and these are the only partitions with corresponding irreducible representations of dimension $< 400$ which
have new identities.
\end{theorem}


\section{A new multilinear identity for representation 51} \label{reptheorysection}

Representation 51 is the smaller of the two representations with new identities in Table \ref{degree11ranks}.
In this section we obtain an explicit form of a new identity for this representation.
(Similar computations could be performed for representation 52.)

From the computations in the previous section we obtain two matrices:
  \begin{itemize}
  \item
  \texttt{oldmat} of size $1020 \times 1056$: this full rank matrix contains the rows representing
  the symmetries of the association types and the consequences of $I(a,b,c,d,e,f,g)$
  for the representation corresponding to partition $\lambda = 2^5 1$.
  \item
  \texttt{allmat} of size $1021 \times 1056$: this full rank matrix contains the rows representing
  all the polynomial identities satisfied by the ternary commutator for the representation corresponding
  to partition $\lambda = 2^5 1$.
  \end{itemize}
The row space of \texttt{oldmat} is a subspace of the row space of \texttt{allmat}.
For a matrix $A$ in row canonical form, we write $\texttt{leading}(A)$ for the set of column indices
for those columns which contain the leading 1 of some row.
We have
  \allowdisplaybreaks
  \begin{align*}
  &
  \texttt{leading}(\texttt{oldmat}) \subset \texttt{leading}(\texttt{allmat}),
  \\
  &
  \texttt{leading}(\texttt{allmat}) - \texttt{leading}(\texttt{oldmat}) = \{ \, 251 \, \}.
  \end{align*}
The row of \texttt{allmat} which has its leading 1 in column 251 is row 246;
this is the row which represents the new identity.
This row has 24 nonzero entries, with 16 distinct integer coefficients:
  \[
  -432, \, -60, \, -36, \, -34, \, -24, \, -9, \, 9, \, 18, \, 24, \, 36, \, 54, \, 72, \, 96, \, 108, \, 144, \, 216.
  \]
The columns of $\texttt{allmat}$ correspond to 8 blocks of length $d_\lambda = 132$; the blocks correspond
to the association types \eqref{types} and the columns in each block correspond to the standard tableaus for partition $\lambda = 2^5 1$
in lexicographical order:
  \[
  \begin{array}{|r|r|}
  \hline
   1 &  2 \\
   \hline
   3 &  4 \\
   \hline
   5 &  6 \\
   \hline
   7 &  8 \\
   \hline
   9 & 10 \\
   \hline
  11 \\
  \cline{1-1}
  \end{array}
  \qquad
  \begin{array}{|r|r|}
  \hline
   1 &  2 \\
   \hline
   3 &  4 \\
   \hline
   5 &  6 \\
   \hline
   7 &  8 \\
   \hline
   9 & 11 \\
   \hline
  10 \\
  \cline{1-1}
  \end{array}
  \qquad
  \begin{array}{|r|r|}
  \hline
   1 &  2 \\
   \hline
   3 &  4 \\
   \hline
   5 &  6 \\
   \hline
   7 &  9 \\
   \hline
   8 & 10 \\
   \hline
  11 \\
  \cline{1-1}
  \end{array}
  \qquad
  \cdots
  \qquad
  \begin{array}{|r|r|}
  \hline
   1 &  6 \\
   \hline
   2 &  8 \\
   \hline
   3 &  9 \\
   \hline
   4 & 10 \\
   \hline
   5 & 11 \\
   \hline
   7 \\
  \cline{1-1}
  \end{array}
  \qquad
  \begin{array}{|r|r|}
  \hline
   1 &  7 \\
   \hline
   2 &  8 \\
   \hline
   3 &  9 \\
   \hline
   4 & 10 \\
   \hline
   5 & 11 \\
   \hline
   6 \\
  \cline{1-1}
  \end{array}
  \]
Table \ref{newrow} gives complete information about the row representing the new identity,
where $t$ is the association type, $j$ is the tableau index, and $c$ is the coefficient;
the standard tableaus are given in flattened form as a sequence of rows.

  \begin{table}
  \[
  \begin{array}{rrrrr|rr|rr|rr|rr|rr}
  \text{column} &\quad t &\quad j &\multicolumn{11}{c}{\quad\text{standard tableau}} &\quad c \\
  \midrule
   251  &\quad  2  &\quad  119  &\quad   1 &  5 &  2 &  6 &  3 &  8 &  4 &  9 &  7 & 11 & 10  &\quad   72 \\
   253  &\quad  2  &\quad  121  &\quad   1 &  5 &  2 &  6 &  3 &  9 &  4 & 10 &  7 & 11 &  8  &\quad  -36 \\
   361  &\quad  3  &\quad   97  &\quad   1 &  4 &  2 &  5 &  3 &  8 &  6 & 10 &  7 & 11 &  9  &\quad   54 \\
   378  &\quad  3  &\quad  114  &\quad   1 &  5 &  2 &  6 &  3 &  7 &  4 &  8 &  9 & 11 & 10  &\quad  144 \\
   388  &\quad  3  &\quad  124  &\quad   1 &  5 &  2 &  7 &  3 &  8 &  4 & 10 &  6 & 11 &  9  &\quad  216 \\
   393  &\quad  3  &\quad  129  &\quad   1 &  6 &  2 &  7 &  3 &  8 &  4 & 10 &  5 & 11 &  9  &\quad  -60 \\
   396  &\quad  3  &\quad  132  &\quad   1 &  7 &  2 &  8 &  3 &  9 &  4 & 10 &  5 & 11 &  6  &\quad   36 \\
   528  &\quad  4  &\quad  132  &\quad   1 &  7 &  2 &  8 &  3 &  9 &  4 & 10 &  5 & 11 &  6  &\quad    9 \\
   622  &\quad  5  &\quad   94  &\quad   1 &  4 &  2 &  5 &  3 &  7 &  6 & 10 &  8 & 11 &  9  &\quad  108 \\
   623  &\quad  5  &\quad   95  &\quad   1 &  4 &  2 &  5 &  3 &  8 &  6 &  9 &  7 & 10 & 11  &\quad  -36 \\
   626  &\quad  5  &\quad   98  &\quad   1 &  4 &  2 &  5 &  3 &  9 &  6 & 10 &  7 & 11 &  8  &\quad  108 \\
   645  &\quad  5  &\quad  117  &\quad   1 &  5 &  2 &  6 &  3 &  7 &  4 & 10 &  8 & 11 &  9  &\quad  216 \\
   653  &\quad  5  &\quad  125  &\quad   1 &  5 &  2 &  7 &  3 &  9 &  4 & 10 &  6 & 11 &  8  &\quad -432 \\
   655  &\quad  5  &\quad  127  &\quad   1 &  6 &  2 &  7 &  3 &  8 &  4 &  9 &  5 & 10 & 11  &\quad   24 \\
   658  &\quad  5  &\quad  130  &\quad   1 &  6 &  2 &  7 &  3 &  9 &  4 & 10 &  5 & 11 &  8  &\quad  144 \\
   660  &\quad  5  &\quad  132  &\quad   1 &  7 &  2 &  8 &  3 &  9 &  4 & 10 &  5 & 11 &  6  &\quad   96 \\
   778  &\quad  6  &\quad  118  &\quad   1 &  5 &  2 &  6 &  3 &  8 &  4 &  9 &  7 & 10 & 11  &\quad  -24 \\
   781  &\quad  6  &\quad  121  &\quad   1 &  5 &  2 &  6 &  3 &  9 &  4 & 10 &  7 & 11 &  8  &\quad   72 \\
   792  &\quad  6  &\quad  132  &\quad   1 &  7 &  2 &  8 &  3 &  9 &  4 & 10 &  5 & 11 &  6  &\quad  -34 \\
   890  &\quad  7  &\quad   98  &\quad   1 &  4 &  2 &  5 &  3 &  9 &  6 & 10 &  7 & 11 &  8  &\quad   -9 \\
   916  &\quad  7  &\quad  124  &\quad   1 &  5 &  2 &  7 &  3 &  8 &  4 & 10 &  6 & 11 &  9  &\quad  216 \\
   918  &\quad  7  &\quad  126  &\quad   1 &  5 &  2 &  8 &  3 &  9 &  4 & 10 &  6 & 11 &  7  &\quad  108 \\
   924  &\quad  7  &\quad  132  &\quad   1 &  7 &  2 &  8 &  3 &  9 &  4 & 10 &  5 & 11 &  6  &\quad  108 \\
  1056  &\quad  8  &\quad  132  &\quad   1 &  7 &  2 &  8 &  3 &  9 &  4 & 10 &  5 & 11 &  6  &\quad   18 \\
  \midrule
  \end{array}
  \]
  \caption{Row 246 of \texttt{allmat} representing the new identity}
  \label{newrow}
  \end{table}

To convert this data into an explicit identity for the ternary commutator, we use the correspondence
between matrix units in the representation matrices and elements of the group algebra
\cite[Remark 2, p.~2004]{BP1}.
We summarize this result in the general case.
Given a partition $\lambda$ of $n$, let $d_\lambda$ be the dimension of the corresponding irreducible
representation of $S_n$.
For $1 \le i, j \le d_\lambda$ we construct the element of the group algebra $\mathbb{Q} S_n$ corresponding
to the matrix unit $E^\lambda_{ij}$ under the isomorphism of $\mathbb{Q} S_n$ with a direct sum of full matrix
algebras.
Let $T_1, \dots, T_{d_\lambda}$ be the standard tableaus for $\lambda$ in lexicographical order.
For each $i = 1, \dots, d_\lambda$ let $R_i$ (respectively $C_i$) be the subgroup of $S_n$ which leaves the rows
(respectively columns) of $T_i$ fixed as sets.
For $i, j = 1, \dots, d_\lambda$ let $s_{ij}$ be the permutation for which $s_{ij} T_i = T_j$.
We define elements $D_{ij} \in S_n$ as follows:
  \[
  D_{ii} = \frac{d_\lambda}{n!} \sum_{\sigma \in R_i} \sum_{\tau \in C_i} \epsilon(\tau) \, \sigma \, \tau,
  \qquad
  D_{ij} = D_{ii} \, s_{ij}^{-1}.
  \]
These elements in general do not satisfy the multiplication formulas for matrix units.
To obtain the matrix units, let $A^\lambda_\pi$ be the matrix defined by Clifton \cite{Clifton}
for the permutation $\pi$.
For the identity permutation $\iota$, the matrix $A^\lambda_\iota$ is not necessarily the identity matrix,
but it is always invertible.
Let $( a_{ij} )$ be the inverse matrix $( A^\lambda_\iota )^{-1}$; then the element of $\mathbb{Q} S_n$
corresponding to the matrix unit $E^\lambda_{ij}$ is
  \[
  E^\lambda_{ij} \; \longleftrightarrow \; \sum_{k=1}^d a_{jk} D_{ik}.
  \]
We then have the required relations $E^\lambda_{ij} E^\lambda_{k\ell} = \delta_{jk} E^\lambda_{i\ell}$.

We now return to our discussion of the new identity in degree 11 for the ternary commutator.
Since we are dealing with a single identity we may assume that $i = 1$: any row of the representation
matrix can be moved to row 1 by left multiplication by an element of the group algebra.
Moreover, we need to consider only those values of $j$ which appear in Table \ref{newrow}:
  \[
  j = 94, \, 95, \, 97, \, 98, \, 114, \, 117, \, 118, \, 119, \, 121, \, 124, \, 125, \, 126, \, 127, \, 129, \, 130, \, 132.
  \]
We compute the matrix $A^\lambda_\iota$ and find that it has the form $I + U$ where $U$ is a strictly
upper triangular matrix with 262 nonzero entries from the set $\{ \pm 1 \}$.
The inverse matrix $( A^\lambda_\iota )^{-1}$ has the form $I + V$ where $V$ is a strictly
upper triangular matrix with 424 nonzero entries from the set $\{ \pm 1, \pm 2 \}$.
For all except one of the values of $j$ listed above, the corresponding row of $( A^\lambda_\iota )^{-1}$
has only one nonzero entry, which is the diagonal entry 1.
Therefore, for all these values except $j = 118$ the matrix unit is $E_{1j} = D_{1j}$;
the exceptional case is
$E_{1,118} \leftrightarrow D_{1,118} - D_{1,126} + D_{1,131}$.
It remains to apply the association types \eqref{types} to the elements of the group algebra.
Let $f \in \mathbb{Q} S_{11}$ be arbitrary, and let $t = 1, \dots, 8$ be one of the association types.
We regard $f$ as a multilinear associative polynomial in the variables $a_1, \dots, a_{11}$.
We convert this into a polynomial in the ternary commutator by applying association type $t$ to every monomial.
The resulting element of the free alternating ternary algebra will be denoted $[f]_t$.
We can now write down the new identity.

\begin{theorem}
The following multilinear polynomial identity in degree 11 is satisfied by the ternary commutator
and does not follow from the symmetries of the association types and the consequences of $I(a,b,c,d,e,f,g)$:
  \allowdisplaybreaks
  \begin{align*}
  &
  \quad\,
    72 \,\, [ D_{1,119} ]_2
   -36 \,\, [ D_{1,121} ]_2
   +54 \,\, [ D_{1, 97} ]_3
  +144 \,\, [ D_{1,114} ]_3
  +216 \,\, [ D_{1,124} ]_3
  \\
  &
   -60 \,\, [ D_{1,129} ]_3
   +36 \,\, [ D_{1,132} ]_3
    +9 \,\, [ D_{1,132} ]_4
  +108 \,\, [ D_{1, 94} ]_5
   -36 \,\, [ D_{1, 95} ]_5
   \\
   &
  +108 \,\, [ D_{1, 98} ]_5
  +216 \,\, [ D_{1,117} ]_5
  -432 \,\, [ D_{1,125} ]_5
   +24 \,\, [ D_{1,127} ]_5
  +144 \,\, [ D_{1,130} ]_5
  \\
  &
   +96 \,\, [ D_{1,132} ]_5
   -24 \,\, \big( \, [ D_{1,118} ]_6 - [ D_{1,126}  ]_6 + [ D_{1,131} ]_6 \, \big)
   +72 \,\, [ D_{1,121} ]_6
   \\
   &
   -34 \,\, [ D_{1,132} ]_6
    -9 \,\, [ D_{1, 98} ]_7
  +216 \,\, [ D_{1,124} ]_7
  +108 \,\, [ D_{1,126} ]_7
  +108 \,\, [ D_{1,132} ]_7
  \\
  &
   +18 \,\, [ D_{1,132} ]_8.
  \end{align*}
This identity implies all the new identities for partition $2^5 1$.
\end{theorem}


\section{A new non-multilinear identity for representation 51}

For partition $\lambda = 2^5 1$ we expect that there will be a new identity in which each term consists of
a permutation of the multiset $a^2 b^2 c^2 d^2 e^2 f$ with one of the eight association types \eqref{types}.
The total number of such permutations is $\binom{11}{2,2,2,2,2,1} = 1247400$.
These permutations will be called ``associative monomials''; they form a basis of the homogeneous subspace
$A_\delta$ of the free
associative algebra on six generators with multidegree $\delta = (2,2,2,2,2,1)$.
The total number of alternating ternary monomials in each association type can be
determined by direct enumeration.  If $p_1 \cdots p_{11}$ denotes an associative monomial, then
the eight association types require the following conditions, where $<$ denotes lexicographical order:
  \begin{equation}
  \label{nonlineartypes} \tag{$\ast\ast$}
  \left\{
  \begin{array}{rl}
  1\colon
  &\;
  p_1 < p_2 < p_3, \quad p_4 < p_5, \quad p_6 < p_7, \quad p_8 < p_9, \quad p_{10} < p_{11}
  \\[3pt]
  2\colon
  &\;
  p_1 < p_2 < p_3, \quad p_4 < p_5 < p_6, \quad p_8 < p_9, \quad p_{10} < p_{11},
  \\
  &\;
  p_1 p_2 p_3 < p_4 p_5 p_6
  \\[3pt]
  3\colon
  &\;
  p_1 < p_2 < p_3, \quad p_4 < p_5, \quad p_6 < p_7 < p_8, \quad p_{10} < p_{11}
  \\[3pt]
  4\colon
  &\;
  p_1 < p_2 < p_3, \quad p_4 < p_5 < p_6, \quad p_7 < p_8 < p_9, \quad p_{10} < p_{11},
  \\
  &\;
  p_1 p_2 p_3 < p_4 p_5 p_6 < p_7 p_8 p_9
  \\[3pt]
  5\colon
  &\;
  p_1 < p_2 < p_3, \quad p_4 < p_5, \quad p_6 < p_7, \quad p_8 < p_9 < p_{10}
  \\[3pt]
  6\colon
  &\;
  p_1 < p_2 < p_3, \quad p_4 < p_5 < p_6, \quad p_8 < p_9 < p_{10},
  \\
  &\;
  p_1 p_2 p_3 < p_4 p_5 p_6
  \\[3pt]
  7\colon
  &\;
  p_1 < p_2 < p_3, \quad p_4 < p_5, \quad p_6 < p_7 < p_8, \quad p_9 < p_{10},
  \\
  &\;
  p_1 p_2 p_3 p_4 p_5 < p_6 p_7 p_8 p_9 p_{10}
  \\[3pt]
  8\colon
  &\;
  p_1 < p_2 < p_3, \quad p_4 < p_5, \quad p_6 < p_7 < p_8, \quad p_9 < p_{10} < p_{11},
  \\
  &\;
  p_6 p_7 p_8 < p_9 p_{10} p_{11}
  \end{array}
  \right.
  \end{equation}
The total number of permutations satisfying these conditions is
  \[
  6720 + 1980 + 4010 + 180 + 4010 + 1190 + 2000 + 550 = 20640.
  \]
The resulting bracketed permutations will be called ``nonassociative monomials'';
they form a basis of the homogeneous subspace $N_\delta$ of the free alternating ternary
algebra on six generators with multidegree $\delta$.

At this point, we would like to construct a matrix of size $1247400 \times 20640$
in which the $(i,j)$ entry is the coefficient of the $i$-th associative monomial
in the expansion of the $j$-th nonassociative monomial; as before, by expansion we mean
repeated application of the alternating ternary sum.
This matrix represents the ``expansion map'' $E_\delta\colon N_\delta \to A_\delta$ with
respect to the bases of associative and nonassociative monomials.
The polynomial identities satisfied by the ternary commutator are the (nonzero) vectors
in the kernel $K_\delta$ of this linear map.
The matrix representing $E_\delta$ is very sparse, since each expansion contains only 7776 terms;
more than 99\% of the entries are 0.
However, processing a matrix of this size is not practical.
We therefore begin by storing the expansions of the nonassociative monomials in a
matrix of size $7776 \times 20640$; the $(i,j)$ entry contains the $i$-th term of
the $j$-th expansion in the form $\pm k$.  The sign $\pm 1$ is the coefficient of
the term and the absolute value $k$ is the lexicographical index of the associative monomial.

We now observe that $1247400 = 77 \cdot 16200$.  We construct a matrix with an
upper  block of size $20640 \times 20640$ and a lower  block of
size $16200 \times 20640$, and initialize it to zero.  We then perform the following
iteration for $\ell = 1, \dots, 77$:
  \begin{itemize}
  \item
  For each column index $j$, extract the terms of the corresponding expansion whose
  indices $k$ lie in the range $16200 (\ell{-}1) < k \le 16200 \ell$.
  \item
  Store the corresponding coefficients in the appropriate row of the lower block;
  index $k$ goes to row $k {-} 16200 (\ell{-}1)$.
  \item
  After all the columns have been processed, and the lower block has been filled,
  compute the row canonical form.  (The lower block is now zero.)
  \end{itemize}
At the end of this iteration, the nullspace of the matrix contains the coefficient
vectors of the polynomial identities satisfied by the ternary commutator.
The rank of the matrix is 19964, and so the nullity is 676.
We compute the canonical basis of the nullspace from the row canonical form by setting
the free variables equal to the standard basis vectors in dimension 676 and solving for the
leading variables.
We summarize this discussion in the following lemma.

\begin{lemma} \label{51nullspacelemma}
The kernel $K_\delta$ of the linear map $E_\delta\colon N_\delta \to A_\delta$ has dimension 676.
\end{lemma}

The next step is to determine the subspace $L_\delta \subset K_\delta$ consisting of
the polynomial identities
which are consequences of the known identity $I(a,b,c,d,e,f,g)$ in degree 7.  (The consequences of
the alternating properties in degree 3 have already been excluded by our choice \eqref{nonlineartypes} of
nonassociative monomials.)
For each consequence of $I(a,b,c,d,e,f,g)$ in Table \ref{degree11liftings}, we must determine
the corresponding substitutions of the variables $a^2 b^2 c^2 d^2 e^2 f$, recalling that $[a,b,c]$
alternates in all three arguments and $I(a,b,c,d,e,f,g)$ alternates in $b,c,d,e,f,g$.
If $q_1 \cdots q_{11}$ denotes an associative monomial, then the eight consequences require the following conditions,
where $<$ denotes lexicographical order:
  \allowdisplaybreaks
 \begin{align*}
  1\colon
  &\quad
  q_1 < q_{10} < q_{11}, \quad q_8 < q_9, \quad q_2 < q_3 < q_4 < q_5 < q_6 < q_7
  \\
  2\colon
  &\quad
  q_1 < q_8 < q_9, \quad q_2 < q_{10} < q_{11}, \quad q_3 < q_4 < q_5 < q_6 < q_7
  \\
  3\colon
  &\quad
  q_1 < q_8 < q_9, \quad q_2 < q_3 < q_4 < q_5 < q_6 < q_7, \quad q_{10} < q_{11}
  \\
  4\colon
  &\quad
  q_2 < q_{10} < q_{11}, \quad q_8 < q_9, \quad q_3 < q_4 < q_5 < q_6 < q_7
  \\
  5\colon
  &\quad
  q_2 < q_8 < q_9, \quad q_3 < q_{10} < q_{11}, \quad q_4 < q_5 < q_6 < q_7,
  \quad
  q_2 q_8 q_9 < q_3 q_{10} q_{11}
  \\
  6\colon
  &\quad
  q_2 < q_8 < q_9, \quad q_3 < q_4 < q_5 < q_6 < q_7, \quad q_{10} < q_{11}
  \\
  7\colon
  &\quad
  q_2 < q_3 < q_4 < q_5 < q_6 < q_7, \quad q_8 < q_{10} < q_{11}
  \\
  8\colon
  &\quad
  q_2 < q_3 < q_4 < q_5 < q_6 < q_7, \quad q_8 < q_9, \quad q_{10} < q_{11}
  \end{align*}
The total number of substitutions satisfying these conditions is
  \[
  10 + 50 + 10 + 170 + 215 + 170 + 20 + 30 = 675.
  \]
We note that this number is exactly one less than the dimension of $K_\delta$.
Each of these consequences of $I(a,b,c,d,e,f,g)$ expands to a linear combination of
terms which consist of a sign and a permutation of $a^2 b^2 c^2 d^2 e^2 f$ with one of the eight
association types \eqref{types}.
We construct a matrix $M$ of size $676 \times 20640$ and fill the first 675 rows with
the coefficient vectors of these consequences of $I(a,b,c,d,e,f,g)$.
We compute the row canonical form, and find that the rank is 675, so these substitutions are
linearly independent.
We summarize this discussion in the following lemma.

\begin{lemma} \label{51liftedlemma}
The subspace $L_\delta \subset K_\delta$ spanned by the consequences of the identity $I(a,b,c,d,e,f,g)$
has dimension 675.
\end{lemma}

It follows that any complementary subspace to $L_\delta$ inside $K_\delta$ has dimension 1;
this agrees with the result for partition $2^5 1$ from Table \ref{degree11ranks}.
Our next task is to find a basis for this complementary subspace.

In each vector of the canonical basis of $K_\delta$, we clear denominators and cancel common factors,
so that the components are relatively prime integers.
We then sort the vectors by increasing Euclidean length;
the minimum square length is 60 and the maximum is 79134357.
In these 676 vectors, the minimum number of nonzero components is 58 and the maximum is 15901;
the minimum number of distinct coefficients is 2 and the maximum is 509.
We copy these row vectors one at a time to the last row of the matrix $M$; after each row, we reduce the matrix.
The vector which increases the rank from 675 to 676 corresponds to item 585 in the sorted list,
which is item 241 of the original (unsorted) canonical basis of the nullspace.
This vector has 10292 nonzero components.

\begin{theorem}
There exists a non-multilinear polynomial identity in degree 11 satisfied by the ternary commutator which is not
a consequence of the identities of lower degree.
This identity has 10292 terms in the six variables $a^2 b^2 c^2 d^2 e^2 f$; it involves only association
types {1, 2, 3, 5, 6} and its coefficients are
$\pm 1$, $\pm 2$, $\pm 3$, $\pm 4$, $\pm 5$, $\pm 6$, $\pm 7$, $\pm 8$, $\pm 9$, $\pm 10$, $\pm 11$, $-12$, $13$.
\end{theorem}

To complete the computational verification of this theorem, we proceed as follows.
We first check by direct expansion that this polynomial is in fact an identity satisfied
by the ternary commutator.  We create a vector of length 1247400, initialized to 0; this vector
represents a linear combination of the 1247400 associative monomials in the variables $a^2 b^2 c^2 d^2 e^2 f$.
We then expand each of the 10292 terms of the polynomial using the alternating ternary sum; for
each term we obtain a linear combination of 7776 terms, where each term is ($\pm$) one of the
associative monomials.  For each term in the expansion, we add the appropriate coefficient to the
corresponding component of the array.  After all terms of the polynomial have been expanded and added
to the vector, every component of the vector is 0.

Finally, we need to verify that this polynomial increases the rank from 1020 to 1021 in the
irreducible representation for partition $2^5 1$; see Table \ref{degree11ranks}.
We first linearize the polynomial: each term $\cdots a \cdots a \cdots$ has two occurrences of the variable $a$ and produces two terms
$\cdots a \cdots g \cdots$ and $\cdots g \cdots a \cdots$; then we similarly replace $bb$ by $bh$ and $hb$,
$cc$ by $ci$ and $ic$, $dd$ by $dj$ and $jd$, $ee$ by $ek$ and $ke$.  Each term of the original
non-multilinear polynomial produces 32 terms in the linearized polynomial, giving a total of $32 \cdot 10292 =
329344$ terms.

We now use the representation theory of the symmetric group as described in Section \ref{reptheorysection}
to redo the computation for partition $2^5 1$ with dimension $d_\lambda = 132$.
We construct a $52 d_\lambda \times 8 d_\lambda$ matrix consisting of $d_\lambda \times d_\lambda$ blocks;
the first 43 rows of blocks contain the representation matrices for the symmetries of the association
types, and the next 8 rows of blocks contain the representation matrices for the consequences of $I(a,b,c,d,e,f,g)$.
The last row of blocks contains the representation matrices for the linearized form of the new identity.
We compute the row canonical form of this matrix and find that its rank is 1021, as required.
Furthermore, we verify that the nonzero rows of this matrix coincide exactly with the $1021 \times 1056$
matrix \texttt{allmat} obtained from the expansions of the association types.


\section*{Acknowledgements}

Murray Bremner was partially supported by a Discovery Grant from NSERC of Canada,
and Luiz Peresi was partially supported by a grant from CNPq of Brazil.


\end{document}